\documentclass[a4,12pt]{article}
\usepackage{bbm}
\usepackage{amsmath}
\usepackage{amssymb}
\usepackage{amsfonts}
\usepackage{cite}
\begin{document}

\noindent

\begin{center}

  {\LARGE \bf  Complete Affine K$\ddot{a}$hler  Manifolds }

  \end{center}

\noindent
\begin{center}

   {  \large Fang Jia and \large An-Min Li  }\footnote{The first author
   is partially supported by NSFC 10871136 grand. The second  author
   is partially supported by NKBRPC(2006 CB805905),NSFC 10631050,RFDP and
   AvH}\\[5pt]

\end{center}
{\bf Abstract.} In this paper we prove that for a complete,
connected and oriented Affine K\"{a}hler manifold $(M,G)$ of
dimension $n,$ if it is affine
  K\"ahler  Ricci flat or if the affine K$\ddot{a}$hler
 scalar curvature $S\equiv0,$ ($n\leq 5$), then the affine K\"{a}ler metric is flat.

 \vskip 0.1in \noindent MSC 2000: 53A15\\Keywords: K\"{a}hler affine  manifold,
 K\"ahler affine Ricci flat,K$\ddot{a}$hler
 affine scalar curvature.
  \vskip 0.1in\noindent{\large \bf Introduction } \vskip
0.1in \noindent

  An affine manifold is a manifold which can be covered by coordinate charts so that the
 coordinate transformations are given by invertible affine transformations. Let $M$ be such
 affine manifold. We shall always assume that our coordinate systems are chosen as above and we
 call it affine coordinates. Let $M$ be an affine manifold.  An affine K$\ddot{a}$hler metric on $M$
 is a Riemannian metric on $M$ such that locally, for affine coordinates $(x_{1},x_{2},\cdots,
 x_{n})$, there is a potential $f$ such that\vskip 0.1in \noindent
$$G_{ij}=\frac{\partial^{2} f}{\partial x_{i}\partial x_{j}}.$$\vskip 0.1in \noindent The pair
 $\left(M, G\right)$ is called an affine K$\ddot{a}$hler manifold.\vskip
 0.1in \noindent

 It is easy to see that the tangent bundle of an affine
  manifold is naturally a complex manifold.
  For each coordinate chart $(x_{1},x_{2},\cdots,x_{n})$,
  if we write a tangent vector of $M$ as
$ \sum y_{i}\frac{\partial}{\partial x_{i}},$ then
$$z_{i}=x_{i}+\sqrt{-1}y_{i},\;i=1,2,\cdots,n$$ are local
holomorphic coordinates of $TM.$ The affine K$\ddot{a}$hler metric
 naturally extends to be a  K$\ddot{a}$hler metric of the complex manifold.
 The Ricci curvature and the scalar curvature  of this K$\ddot{a}$hler metric
 are given respectively by\vskip 0.1in \noindent
 $$K_{ij}=-\sum \frac{\partial^{2}}{\partial x_{i}\partial x_{j}}\left(\log \det\left(f_{kl}\right)\right),$$
$$S=-\frac{1}{2}\sum^{n}_{i=1}\sum^{n}_{j=1}f^{ij}\frac{\partial^{2}\log \det\left(
  f_{kl}\right) }{\partial x_{i}\partial x_{j}}
 ,$$where $f_{kl}=\frac{\partial^{2}f}{\partial x_{k}\partial x_{l}}.$
Following Cheng and Yau ([2]) we call $K_{ij}$ and $S$ the affine
K$\ddot{a}$hler Ricci curvature and
 the affine K$\ddot{a}$hler scalar curvature of $(M,G)$ respectively. We say that the affine K$\ddot{a}$hler
  metric $G$ is Einstein if its Ricci tensor is a scalar multiple of the affine K$\ddot{a}$hler
  metric, that is \vskip 0.1in \noindent
$$- \frac{\partial^{2}}{\partial x_{i}\partial x_{j}}\left(\log \det\left(f_{kl}\right)\right)
 = a f_{ij} ,$$\vskip 0.1in \noindent where
$a$ is
  a constant. In particular, if $a=0$ then we call $(M,G)$
affine  K\"ahler Ricci flat.

\vskip 0.1in \noindent
 Our main results can be stated as follows:
 \vskip 0.1in \noindent{\bf Theorem 1.} {\it
  Let $(M,G)$ be an n-dimensional complete, connected and oriented $C^{\infty}$ affine K\"ahler Ricci flat
   manifold.  Then the affine K$\ddot{a}$hler metric is flat.}\vskip 0.1in \noindent
 {\bf Theorem 2.} {\it Let $(M,G)$ be a complete, connected and oriented $C^{\infty}$
affine K$\ddot{a}$hler manifold of dimension $n$.   If the affine K$\ddot{a}$hler
 scalar curvature $S\equiv0,$ then , for $n\leq 5$, the affine K$\ddot{a}$hler metric is flat.}
  \vskip 0.1in \noindent
  As consequences we have
   \vskip 0.1in \noindent{\bf Theorem 3.} {\it
  Let $f(x_1,...,x_n)$ be a smooth and strictly convex function defined in $\Omega \subset R^n$. If the affine K\"ahler Ricci curvature is identically $0$, and if the Calabi metric is complete, then $f$ must be
  a quadratic polynomial.
  } \vskip 0.1in \noindent
 {\bf Theorem 4.} {\it Let $f(x_1,...,x_n)$ be a smooth and strictly convex function defined in $\Omega \subset R^n$. If the affine K$\ddot{a}$hler
 scalar curvature $S\equiv0,$ and if the Calabi metric is complete, then , for $n\leq 5$, $f$ must be
  a quadratic polynomial.
  }\vskip 0.1in \noindent

{\bf Remark.} This is an unpublished paper that was finished in
2005. Since then, the formula of $\Delta \Phi$ (cf. Proposition
1) has been used frequently in various circumstances. Since it is
used in our recent papers again (see \cite{9,10}),   we decide to put this orginal version (with
slight revision) on Arxiv. 

\section{Fundamental formulas }\vskip 0.1in \noindent
Let $(M,G)$ be an affine K\"ahler  manifold. Choose a local affine
coordinate system $(x_1,...,x_n).$ Let $f(x)$ be a local potential
function of $G.$ Then $f$ is locally strictly convex function and
$$G=\sum_{i,j} f_{ij}d x_i dx_j .$$
 We recall some fundamental
facts on the Riemannian manifold $(M, G)$ (cf. \cite{9}). The
Levi-Civita connection is given by
$$\Gamma_{ij}^k=\frac{1}{2}\sum {f^{kl}}{f_{ijl}}.$$
The Fubini-Pick tensor is
$$A_{ijk}=- \frac{1}{2}f_{ijk}.$$
Then the  curvature tensor and the Ricci tensor are
\begin{eqnarray}
R_{ijkl} &=& \sum f^{mh}(A_{jkm}A_{hil}-A_{ikm}A_{hjl})
\nonumber\\\label{eqn_0}
R_{ik}&=&\sum f^{mh}f^{jl}(A_{jkm}A_{hil}-A_{ikm}A_{hjl}).
\end{eqnarray}\vskip 0.1in \noindent
Let $ \rho=\left[\det(f_{ij})\right]^{-\frac{1}{n+2}}. $ Set
\begin{equation}\label{eqn_2.1}
\Phi=\frac{\|\nabla\rho\|^2_G}{\rho^2} \end{equation}
\begin{equation}\label{eqn_2.2}
4n(n-1)J=\sum f^{il}f^{jm}f^{kn}f_{ijk}f_{lmn}.
\end{equation}
It is easy to check that $\Phi$ and $J$ are independent of the
choice of the affine coordinate systems. Hence they are invariants
globally defined on $M.$ If $\Phi\equiv0$ then $\rho = constant.$ It is well known that
(see \cite{6}) \begin{equation}\label{eqn_3}\triangle J\geq 2(n+1)J^2.\end{equation} Here and later
the Laplacian and the covariant
  differentiation with respect to the  metric $G$ will be denoted by``$\Delta $ "
  and ``," respectively.

\section{\bf Estimate for $\triangle  \Phi$}\vskip 0.1in \noindent

 In this section we calculate $\triangle  \Phi$ for  affine K$\ddot{a}$hler manifold with $S=0$ and
 affine  K$\ddot{a}$hler Ricci flat manifold.\vskip 0.1in\noindent {\bf Proposition 1.}
 Let $(M,G)$ be an n-dimensional ,
 connected and oriented $C^{\infty}$ affine K$\ddot{a}$hler Ricci flat manifold.   Then the following estimate holds
$$\triangle \Phi\geq\frac{n}{n-1}\sum\frac{||\nabla\Phi||_G^2}{\Phi}+
  \frac{n^{2}-3n-10}{2(n-1)}<\nabla \Phi,\nabla \log \rho>_G +\frac{(n+2)^{2}}{n-1}\Phi^{2}.$$
\vskip 0.1in\noindent {\bf Proof.}
 Let $p\in M$ be
any fixed point. Choose an affine coordinate neighborhood $
\left\{U,\varphi\right\}$ with
 $p\in U$.
  We have:\vskip 0.1in \noindent
$$-\frac{\partial^{2}}{\partial x_{i}\partial x_{j}}\left(\log \det\left(f_{kl}\right)\right)=
  (n+2)\left(\frac{\rho_{ij}}{\rho}-\frac{\rho_{i}}{\rho}\frac{\rho_{j}}{\rho}\right),$$
  \vskip 0.1in \noindent where $\rho_{i}=\frac{\partial\rho}{\partial x_{i}}$ and $\rho_{ij}=
  \frac{\partial^{2}\rho}{\partial x_{i}\partial x_{j}}.$ Noting
  that
$-\frac{\partial^{2}}{\partial x_{i}\partial x_{j}}\left(\log \det\left(f_{kl}\right)\right)=
  0,$  we obtain \vskip 0.1in \noindent
$$\frac{1}{\rho}\sum f^{ij}\rho_{ij}+\frac{n}{\rho^{2}}\sum f^{ij}\rho_{i}\rho_{j}-\frac{n+1}
  {\rho^{2}}\sum f^{ij}\rho_{i}\rho_{j}=0,$$\vskip 0.1in\noindent where  the matrix $\left( f^{ij}\right)$ denotes
  the inverse matrix of the matrix $\left( f_{ij}\right).$  Then we have
\begin{equation}\label{eqn_4}\triangle
\rho=\frac{n+4}{2}\frac{\left\|\nabla\rho\right\|^{2}_{G}}{\rho}.\end{equation} We choose a local
orthonormal
  frame field of the metric $G$ on $U$. Then
$$\Phi=\sum\frac{(\rho_{,j})^2}{\rho^{2}}, \;\;\;\; \Phi_{,i}=2\sum\frac{\rho_{,j}\rho_{,ji}}{\rho^{2}}-2\rho_{,i}\sum\frac{(\rho_{,j})^2}
  {\rho^{3}},$$
$$\triangle \Phi=2\sum\frac{(\rho_{,ji})^2}{\rho^{2}}+2\sum\frac{\rho_{,j}
  \rho_{,jii}}{\rho^{2}}-8\sum\frac{\rho_{,j}\rho_{,i}\rho_{,ji}}{\rho^{3}}-(n-2)\frac{\left(
  \sum\rho^{2}_{,j}\right)^{2}}{\rho^{4}},$$ where we used \eqref{eqn_4}. In the case $\Phi(p)=0$,
  it is easy to get , at $p,$
\begin{equation}\label{eqn_5}\triangle \Phi\geq2\sum\frac{(\rho_{,ij})^2}{\rho^{2}}.\end{equation}  Now we assume
  that $\Phi (p)\neq 0$. Choose a local orthonormal frame field of the
  metric$g$ on $U$ such that $$\rho_{,1}(p)=\left\|\nabla\rho\right\|_{G}(p)>0,\;\;\;\rho_{,i}(p)
  =0 \;\;\;\; for\;\; all\;\; i>1.$$ Then
\begin{equation}\label{eqn_6}\triangle \Phi=2\sum\frac{(\rho_{,ij})^2}{\rho^{2}}+2\sum\frac{\rho_{,j}\rho_{,jii}}
  {\rho^{2}}-8\frac{(\rho_{,1})^2\rho_{,11}}{\rho^{3}}-(n-2)\frac{(\rho_{,1})^4}{\rho^{4}}.
  \end{equation}
Applying an elementary
inequality
$$a_1^2+a_2^2+\dots+a_{n-1}^2\geq\frac{(a_1+a_2+\dots+a_{n-1})^2}{n-1}$$
and \eqref{eqn_4}, we obtain
$$2\frac{\sum (\rho_{,ij})^2}{\rho^{2}}\geq
2\frac{(\rho_{,11})^2}{\rho^{2}}+4
\frac{\sum_{i>1}(\rho_{,1i})^2}{\rho^{2}}+2
\frac{\sum_{i>1}(\rho_{,ii})^2}{\rho^{2}} $$$$\geq
2\frac{(\rho_{,11})^2}{\rho^{2}}+
4\frac{\sum_{i>1}(\rho_{,1i})^2}{\rho^{2}} +\frac{2}{n-1}
\frac{(\Delta\rho -\rho_{,11})^2}{\rho^{2}} $$
\begin{equation}\label{eqn_7}
\geq\frac{2n}{n-1}\frac{(\rho_{,11})^2}
{\rho^{2}}+4\frac{\sum_{i>1}(\rho_{,1i})^2}{\rho^{2}}-2\frac{n+4}{n-1}
\frac{(\rho_{,1})^2\rho_{,11}}{\rho^{3}}+\frac{(n+4)^{2}}{2(n-1)}
\frac{(\rho_{,1})^4}{\rho^{4}}.\end{equation}
An application of the Ricci identity shows that\vskip 0.1in \noindent
\begin{equation}\label{eqn_8}\frac{2}{\rho^{2}}\sum\rho_{,j}\rho_{,jii}=2(n+4)\frac{(\rho_{,1})^2\rho_{,11}}{\rho^{3}}
  -(n+4)\frac{(\rho_{,1})^4}{\rho^{4}}+2R_{11}\frac{(\rho_{,1})^2}{\rho^{2}}.\end{equation}
  Substituting \eqref{eqn_7} and \eqref{eqn_8} into \eqref{eqn_6} we obtain

\begin{equation}\label{eqn_9}\triangle \Phi\geq\frac{2n}{n-1}\sum\frac{(\rho_{,11})^2}{\rho^{2}}
  +\left(2n-2\frac{n+4}{n-1}\right)\frac{(\rho_{,1})^2\rho_{,11}}{\rho^{3}}\end{equation}
$$+2R_{11}\frac{(\rho_{,1})^2}{\rho^{2}}+\left(\frac{(n+4)^{2}}{2(n-1)}-2(n+1)\right)
  \frac{(\rho_{,1})^4}{\rho^{4}}+4\sum_{i>1}\frac{(\rho_{,1i})^2}{\rho^{2}}.$$ Note that
\begin{equation}\label{eqn_10}\sum\frac{(\Phi_{,i})^2}{\Phi}=4\sum\frac{(\rho_{,1i})^2}{\rho^{2}}-8\frac{(\rho_{,1})^2
  \rho_{,11}}{\rho^{3}}+4\frac{(\rho_{,1})^4}{\rho^{4}},\end{equation}
  Then \eqref{eqn_9} and \eqref{eqn_10} together give us
\begin{eqnarray}\triangle \Phi&\geq&\frac{n}{2(n-1)}\sum\frac{(\Phi_{,i})^2}{\Phi}+\left(\frac{2n-8}{n-1}+2n\right)
  \frac{(\rho_{,1})^2\rho_{,11}}{\rho^{3}}\nonumber\\\label{eqn_11}
&&+2R
_{11}\frac{(\rho_{,1})^2}{\rho^{2}}+\left(\frac{(n+4)^{2}}{2(n-1)}-2(n+1)-\frac{2n}{n-1}
  \right)\frac{(\rho_{,1})^4}{\rho^{4}}.\end{eqnarray}

From $\frac{\partial^2   }{\partial x_i \partial x_j }\log \det(f_{kl})=0$ we easily obtain
$$\rho_{,ij}=\rho_{ij} +A_{ij1}\rho_{,1}=
\frac{\rho_{,i}\rho_{,j}}{\rho} +A_{ij1}\rho_{,1}.$$ Thus we get
\begin{equation}\label{eqn_12} \Phi_{,i}=\frac{2\rho_{,1}\rho_{,1i}}{\rho^2}
-2\frac{\rho_{,i}(\rho_{,1})^2}{\rho^3}=
2A_{i11}\frac{(\rho_{,1})^2}{\rho^2} ,\end{equation} \begin{equation}\label{eqn_13} \quad \frac{\sum
(\Phi_{,i})^2}{\Phi} = 4\sum
(A_{i11})^2\frac{(\rho_{,1})^2}{\rho^2},\;\;\;\; \sum
\Phi_{,i}\frac{\rho_{,i}}{\rho} =
2A_{111}\frac{(\rho_{,1})^3}{\rho^3}.\end{equation}
   By the same argument of \eqref{eqn_7} we have
$$ \sum (f_{ml1})^2\geq (f_{111})^2+2\sum_{i>1}(f_{i11})^2+\sum_{i>1} (f_{ii1})^2
$$$$\geq (f_{111})^2+2\sum_{i>1}(f_{i11})^2
+\frac{1}{n-1}\left(\sum f_{ii1}-f_{111}\right)^2$$$$
\geq\frac{n}{n-1}\sum (f_{i11})^2- \frac{2}{n-1}f_{111} \sum
f_{ii1} + \frac{1}{n-1}\left(\sum f_{ii1}\right)^2 $$
\begin{equation}\label{eqn_14}=\frac{n}{n-1}\sum (f_{i11})^2+\frac{2(n+2)}{n-1}f_{111}\frac{\rho_{1}}{\rho}+
  \frac{(n+2)^{2}}{n-1}\frac{(\rho_{1})^2}{\rho^{2}}.\end{equation} Combing \eqref{eqn_0},
  \eqref{eqn_13} and \eqref{eqn_14}   we have
\begin{equation}\label{eqn_15}2R _{11}(p)\frac{(\rho_{,1})^2}{\rho^{2}}\geq\frac{n}{2(n-1)}\sum\frac{(\Phi_{,i})^2}{\Phi}
  -\frac{(n+2)(n+1)}{2(n-1)}\sum\Phi_{,i}\frac{\rho_{,i}}{\rho}+\frac{(n+2)^{2}}{2(n-1)}
  \frac{(\rho_{,1})^4}{\rho^{4}}.\end{equation} Then
\begin{equation}\label{eqn_16}\triangle \Phi\geq\frac{n}{n-1}\sum\frac{(\Phi_{,i})^2}{\Phi}+
  \frac{n^{2}-3n-10}{2(n-1)}\sum \Phi_{,i}\frac{\rho_{,i}}{\rho}+\frac{(n+2)^{2}}{n-1}\Phi^{2}.
  \end{equation} \vskip 0.1in\noindent\vskip 0.1in\noindent{\bf
  Proposition 2.} Let $(M,G)$ be an n-dimensional, connected and oriented
$C^{\infty}$ affine K$\ddot{a}$hler
manifold with $S\equiv 0.$ We have
\begin{equation}\label{eqn_17}\triangle \Phi
\geq\frac{n}{2(n-1)}\frac{||\nabla\Phi||_G^{2}}{\Phi}+\frac{n^{2}-4}{n-1}
  <\nabla\Phi,\nabla \log \rho>_G +\frac{(n+2)^{2}}{2}\left(\frac{1}{n-1}-\frac{n-1}{4n}\right)
  \Phi^2.\end{equation}
{\bf Proof.} From the proof of Proposition 1 we see that the
equality \eqref{eqn_11} remains hold. On the other hand
$$2R _{11}(p)\frac{(\rho_{,1})^2}{\rho^{2}}=\frac{1}{2}\sum (f_{kj1})^2\frac{(\rho_{,1})^2}
  {\rho^{2}}+\frac{n+2}{2}f_{111}\frac{(\rho_{,1})^3}{\rho^{3}}$$
  $$\geq \frac{1}{2}\left[f^2_{111}+
\frac{1}{n-1}\left(f_{111}+(n+2)\frac{\rho_1}{\rho}\right)^2\right]\frac{\rho^{2}_{,1}}{\rho^{2}}
+\frac{n+2}{2}f_{111}\frac{\rho^{3}_{,1}}{\rho^{3}}
$$
  $$\geq \frac{(n+2)^2}{2(n-1)}\Phi^2 -\frac{(n+2)^{2}(n+1)^2}{8n(n-1)}\Phi^2\geq-\frac{(n+2)^{2}(n-1)}{8n}\frac{(\rho_{1})^4}{\rho^{4}}.$$\vskip 0.1in\noindent
  This combined with \eqref{eqn_11} yields
\begin{eqnarray*}\triangle
\Phi&\geq&\frac{n}{2(n-1)}\sum\frac{(\Phi_{,i})^2}{\Phi}+\left(\frac{2n-8}{n-1}+2n\right)
  \frac{(\rho_{,1})^2\rho_{,11}}{\rho^{3}}
\\&&+2R
_{11}\frac{(\rho_{,1})^2}{\rho^{2}}+\left(\frac{(n+4)^{2}}{2(n-1)}-2(n+1)-
  \frac{2n}{n-1}\right)\frac{(\rho_{,1})^4}{\rho^{4}}\\&\geq&\frac{n}{2(n-1)}\sum\frac{(\Phi_{,i})^2}{\Phi}+\frac{n^{2}-4}{n-1}\sum
  \Phi_{,i}\frac{\rho_{,i}}{\rho}+\frac{(n+2)^{2}}{2}\left(\frac{1}{n-1}-\frac{n-1}{4n}\right)
  \frac{(\rho_{,1})^4}{\rho^{4}}. \end{eqnarray*}\vskip 0.1in \noindent

\section{\bf Proof of Theorems}\vskip 0.1in \noindent

It is well known that an affine complete, parabolic affine hypersphere must be a quadratic
  polynomial. Using \eqref{eqn_3} and the same argument we can get
\vskip 0.1in \noindent {\bf
Lemma 1.} {\it Let $(M,G)$ be a complete, connected and oriented
$C^{\infty}$ affine
 K\"{a}hler manifold of dimension $n$. If $\Phi\equiv0$, then any local potential function $f$ of $G$  must be a quadratic
  polynomial.
  }

  \vskip 0.1in \noindent
{\bf Proof of Theorem 1.} By Lemma 1 it suffices to prove that $\Phi\equiv0$.   Consider the
  function
$$F=(a^{2}-r^{2})^{2}\Phi$$  defined on $B_{a}(p_{0})$. Obviously, $F$
  attains its supremum at some interior point $p^{\ast}$ of $B_{a}(p_{0})$ .     Then
  , at $p^{\ast}$,
\begin{eqnarray}\label{eqn_18}\frac{\Phi_{, i}}{\Phi}-2\frac{(r^{2})_{,
i}}{a^{2}-r^{2}}=0.\\\label{eqn_19}
\frac{\triangle  \Phi}{\Phi}-\sum\frac{(\Phi_{,
i})^2}{\Phi^{2}}-2\sum\frac{\left(r^{2}
  \right)^{2}_{,i}}{(a^{2}-r^{2})^{2}}-2\frac{\triangle  \left(r^{2}\right)}
  {a^{2}-r^{2}}\leq 0.\end{eqnarray} Inserting \eqref{eqn_18} into \eqref{eqn_19}   we get
\begin{equation}\label{eqn_20}\frac{\triangle \Phi}{\Phi}\leq 24\frac{r^{2}}{(a^{2}-r^{2})^{2}}+\frac{4}{a^{2}-
  r^{2}}+4\frac{r\triangle r}{a^{2}-r^{2}}.\end{equation} Denote by $a^{\ast}=
  r(p_{0},p^{\ast})$. In the case $p^{\ast}\neq p_{0}$ we have $a^{\ast}>0$. Let
$$B_{a^{\ast}}(p_{0})=\left\{p\in M|r(p_{0},p)\leq a^{\ast}\right\}.$$
By \eqref{eqn_16} we have
$$\max\limits_{p\in B_{a^{\ast}(p_{0})}}\Phi(p)=\max\limits_{p\in\partial B_{a^{\ast}(p_{0})}}
  \Phi (p).$$\vskip 0.1in \noindent On the other hand, we have $a^{2}-r^{2}=a^{2}-a^{\ast 2}$
   on $\partial B_{a^{\ast}}(p_{0})$, it follows that\vskip 0.1in \noindent
$$\max\limits_{p\in B_{a^{\ast}(p_{0})}}\Phi (p)=\Phi (p^{\ast}).$$\vskip 0.1in \noindent
  Let $ p\in B_{a^{\ast}}(p_{0})$ be any point.
  Then from the definition of $R_{ik}$, we get\vskip 0.1in \noindent
\begin{eqnarray*}R_{ii}(p)&=& \frac{1}{4}\sum f^{jl}f^{hm}(f_{hil}f_{mji}
-f_{hii}f_{mjl})\\&\geq&-\frac{(n+2)^2}{16}\Phi (p)\geq-
\frac{(n+2)^2}{16}\Phi(p^{\ast}).\end{eqnarray*}
  Thus,by Laplacian comparison theorem (see [6] Appendix 2), we obtain\begin{equation}\label{eqn_21}
  r\triangle r\leq (n-1)\left(1+\frac{n+2}{4}\sqrt{\Phi (p^{\ast})}\cdot r\right).\end{equation}
  In the case $p^{\ast}=p_{0}$, we have $r(p_{0},p^{\ast})=0.$
  Consequently, from \eqref{eqn_20} and \eqref{eqn_21}, it follows that
\begin{equation}\label{eqn_22}\frac{\triangle
\Phi}{\Phi}\leq\left(24+\frac{(n-1)^{2}(n+2)^{2}}{4}\right)
  \frac{r^{2}}{(a^{2}-r^{2})^{2}}+\frac{4n}{a^{2}-r^{2}}+\Phi.\end{equation} On the other hand,
  by \eqref{eqn_16} we have
\begin{equation}\label{eqn_23}\frac{\triangle \Phi}{\Phi}\geq-\frac{\left(n^{2}-3n-10\right)^{2}}{(n-1)^{2}}\frac{a^{2}}
  {\left(a^{2}-r^{2}\right)^{2}}+\left(\frac{(n+2)^{2}}{n-1}-1\right)\Phi.\end{equation}
  where we used \eqref{eqn_18}. Inserting \eqref{eqn_23} into \eqref{eqn_22} we get
$$(a^{2}-r^{2})^{2}\Phi\leq C_1(n)a^{2},$$ where $C_1(n)$ is a constant depending only on $n.$ Hence,at any interior point of $B_{\frac{a}{2}}(p_{0}),$ we have
$$\Phi\leq\frac{ 16C_1(n)}{9a^{2}}.$$ \vskip 0.1in\noindent Let
  $a\rightarrow\infty$, then
 $\Phi\equiv 0. $  We complete the proof of Theorem 1.
   \;\;\;\;\;\;\;\;\;$\blacksquare$\vskip 0.1in \noindent

Applying a similar method and using the differential inequality
\eqref{eqn_17} we can prove theorem 2.

\vskip 0.1in \noindent
An-Min Li       \hfill                              Fang Jia \\
Department of Mathematics   \hfill                  Department of Mathematics \\
Sichuan University           \hfill                       Sichuan University \\
Chengdu, Sichuan              \hfill                             Chengdu, Sichuan \\
P.R.CHina                  \hfill                   P.R.China \\
e-mail:math-li@yahoo.com.cn\hfill e-mail:jiafangscu@yahoo.com.cn
 \vskip 0.1in \noindent

\end{document}